\documentclass[final]{siamltex}
\usepackage{graphicx}
\usepackage{amssymb}

\def\Ocal{{\cal O}}
\def\Pcal{{\cal P}}
\def\Dcal{{\cal D}}
\newcommand{\R}{\mbox{I}\!\mbox{R}}
\newcommand{\Z}{{\mathbb Z}}
\newcommand{\N}{{\mathbb N}}
\newcommand{\E}{\mbox{I}\!\mbox{E}\,}
\newcommand{\D}{\mbox{I}\!\mbox{D}}

\title{Non-negativity preserving numerical algorithms
for stochastic differential equations}
\author{
Esteban Moro, Henri Schurz\thanks{Grupo Interdisciplinar de Sistemas Complejos (GISC) and
Departamento de Matem{\'a}ticas, Universidad Carlos III de Madrid,
Avenida de la Universidad 30, 28911 Legan{\'e}s, Spain. (Email {\tt
emoro@math.uc3m.es}), Department of Mathematics,
Southern Illinois University, 1245 Lincoln Drive, Carbondale, IL 62901, USA
and Department of Mathematics and Statistics, Texas Tech University,
Lubbock, TX 79409-1042, USA.
({\tt hschurz@math.siu.edu}). This version is from \today.}
}

\begin{document}

\maketitle

\begin{abstract}
Construction of splitting-step methods and properties of related
non-negativity and boundary preserving numerical algorithms for solving
stochastic differential equations (SDEs) of It{\^o}-type are discussed.
We present convergence proofs for a newly designed splitting-step
algorithm and simulation studies for numerous numerical examples
ranging from stochastic dynamics occurring in asset pricing theory
in mathematical finance (SDEs of CIR and CEV models) to measure-valued
diffusion and superBrownian motion (SPDEs) as met in biology and physics.
\end{abstract}

\begin{keywords}
stochastic differential equations, absorbing boundary, numerical
methods, super-Brownian motion, splitting-step algorithm,
convergence, , simulation
\end{keywords}

\begin{AMS}
65C30, 60H10, 60H15
% 65C30 Probabilistic methods, simulation and stochastic differential equations
%       Stochastic differential and integral equations
% 60H10 Stochastic differential equations
% 60H15 Stochastic analysis
%       Stochastic partial differential equations
\end{AMS}

\pagestyle{myheadings} \thispagestyle{plain} \markboth{E. MORO \&
H. SCHURZ}{NON-NEGATIVITY PRESERVING NUM. ALGORITHMS FOR SDEs}

\section{Introduction and examples}
Stochastic differential equations (SDEs) are fundamental to
describe and understand random phenomena in different areas of
physics, engineering, economics, etc. Particularly, they serve as
model for price fluctuations as in the famous Black-Scholes option
pricing model, description of erratic movements of particles as in
the Langevin equation or spatial processes like the superBrownian
motion. In most cases, explicit solutions of SDEs are very
difficult to obtain and numerical approximations need to be
exploited. In fact a wealth of methods for integrating numerically
SDEs are known and tested
\cite{sdenumeric0,sdenumeric1,sdenumeric2,
sdenumeric3,sdenumeric4}.

In some practical applications autonomous It\^o-type SDEs of the form
\begin{equation}\label{inicio}
dX(t) = f(X(t)) dt + \sigma(X(t)) dW(t),\quad X(0) = X_0 \in \Dcal
\end{equation}
are not well defined unless a boundary condition is additionally
given at the boundaries of the domain $\Dcal$ in which $X(t)$
lives for all $t$ (almost surely). For example, if $X(t)$ is the
price of a stock and (\ref{inicio}) gives its time evolution, then
$\Dcal = [0,\infty)$, where $X(t) = 0$ implies an absorbing or
reflecting state. Or, if $X(t)$ models the number of certain
species in a noisy environment, then $\Dcal = [0,K]$, where $K$ is
an attracting carrying capacity of the environment.

As in most situations, boundary conditions are not needed to state
the related problems (\ref{inicio}) when the boundaries are
unattainable in finite time. This is the case of {\em natural
boundaries}, according to Feller's classification of
one-dimensional diffusions \cite{feller0,feller}. The standard
(unrestricted) Brownian motion on $\R$ is the most obvious example
of diffusion with natural boundaries at infinity. The situation is
different when the solution of (\ref{inicio}) attains the
boundaries in finite time. For example, the Brownian motion on
$\R^+$ in which a boundary condition at $x=0$ needs to be
specified to completely define the solution of it. Typical
boundary conditions in this case are absorbing or reflecting ones
and the solution of (\ref{inicio}) depends on its specific choice, which
is usually taken according to the nature of the problem under
consideration.

This problem of how to handle the boundary conditions also appears
in the numerical approximations of (\ref{inicio}), where the
naturally inherited boundary conditions do need to be incorporated
in the construction of numerical approximations. In fact, standard
numerical methods such as Euler methods may fail to meet the
boundary conditions, see \cite{NonnegativeBims}. This is also true
for higher order converging Taylor-type methods
\begin{equation}\label{generalTaylor}
Y_{n+1} = \sum_{\alpha \in {\mathcal H}} f_\alpha (Y_n) I_\alpha ,
\end{equation}
where $\alpha$ is a multiple index, ${\mathcal H}$ a hierachical
set of multiple indices, $f_\alpha$ coefficient functions and
$I_\alpha$ iterated multiple stochastic integrals as derived from
stochastic Taylor-type expansions (for its origin, see
\cite{wagner-platen78}) along time-discretizations
$$ 0 = t_0 < t_1 < t_2 < ... < t_{n_T} = T $$
with step sizes $\Delta_n t = t_{n+1} - t_n$. For more details,
see \cite{sdenumeric0}, \cite{sdenumeric-1}, \cite{sdenumeric3}.

Numerical time-discretizations of the type (\ref{generalTaylor})
face two kinds of different (though related) problems when a
boundary condition is specified. The first one is to restrict the
values of the numerical approximations to live within the domain
${\cal D}$; secondly, they also have to preserve the character of
the boundary (natural, reflecting, absorbing, etc.) in the
numerical approximations of $X(t)$. To exemplify the numerical
approximation problems, let us consider the well-known square-root
diffusion model of Cox-Ingersoll-Ross \cite{cox}
\begin{equation}\label{coxmodel}
dX(t) = [a + b X(t)] dt + \sigma \sqrt{X(t)} dW(t), \quad X(0) = x_0 \geq 0
\end{equation}
with real parameters $a\geq 0$, $b \in \R$ and $\sigma > 0$, which
is widely used for interest rate modelling or as an alternative to
geometric Brownian motion occurring in the Black-Scholes model of
dynamic asset pricing in mathematical finance. It is well known
that the (strong) solution of (\ref{coxmodel}) is unique and
preserves the non-negativity of the initial data. This property
can be easily implemented in numerical approximations of
(\ref{coxmodel}) by means of balanced implicit methods (see
\cite{mps,sdenumeric2,sdenumeric3}), balanced Milstein methods
\cite{kahl&schurz}, composition methods based on Lie algebra
techniques \cite{misawa}, but fail to incorporate the right
boundary properties in the numerical paths. Other methods, while
converging in the limit $\Delta t \to 0$, do not even conserve the
non-negativity of the solution like the straightforward fixing by
taking $\sqrt{|X(t)|}$ instead of the last term in
(\ref{coxmodel}) (as done in \cite{higham-mao}).

The aim of this paper is to present an alternative strategy based
on {\bf splitting-step methods} to integrate numerically SDEs of
the type (\ref{inicio}) subject to boundary conditions, that both
guarantee that the numerical solutions live in the domain ${\cal
D}$ and that the character of the boundary is preserved (a kind of
numerically compact support property, i.e. the numerical
approximation has to incorporate the fact that, for any positive
$X(t)$, there is a non-zero probability such that the stochastic process
becomes zero at the next time-step). Moreover, the consistency of
this new method shall be mathematically justified by two
convergence theorems, namely one for sufficiently smooth
Lipschitzian coefficients of involved SDEs (Theorem \ref{th:prop})
and another, more complicated one for non-Lipschitzian SDEs
(Theorem \ref{th:sprop}). By those theorems we are going to
establish the same convergence rates as standard, so far known
methods possess. The second Theorem \ref{th:sprop} covers the
square-root case as exhibited by (\ref{coxmodel}) too, and we
present an alternative boundary- and positivity-preserving
numerical splitting algorithm and its proof to that in
\cite{higham-mao} without using substitutions such as
$\sqrt{|X(t)|}$ in diffusion terms of model (\ref{coxmodel}).

The paper is structured as follows. After this introduction we
introduce and discuss the fairly general splitting-step algorithm
in Section 2. Section 3 reports on general convergence theorems
and its mathematical proof. Several numerical experiments in
Section 4 - 6 strongly support the suggested splitting algorithm and
results from previous sections by ordinary SDEs. In Section 4 we
present some simple models and related experiments. Section 5
is devoted to the simulation of processes as often met in dynamic
asset pricing in mathematical finance. In Section 6 we mention
some applications to measure-valued diffusions and
reaction-diffusion equations demonstrating the potential range of
the splitting-step algorithm to numerical treatment of random
PDEs as well. Section 7 concludes with some supplemental
remarks. Eventually, there is a small appendix on aspects of random
number generation.

\section{Splitting-step algorithm}
The general structure of a splitting-step algorithm, which is
based on the idea in \cite{pechenik} is as follows. Suppose that
the more general equation which is to be integrated is of the form
\begin{equation}\label{th11-1}
d X(t) = [\alpha(X(t),t) + \beta(X(t),t)] dt + \sigma(X(t),t) dW(t) .
\end{equation}
We then decompose the above equation into the two equations
\begin{eqnarray}
d X_1(t)&=&\beta(X_1(t),t) dt + \sigma(X_1(t),t) dW(t) , \label{ss-1} \\
d X_2(t)&=&\alpha(X_2(t),t) dt \label{ss-2} ,
\end{eqnarray}
where the splitting is done assuming that we know the exact
strong solution for $X_1(t)$ or the conditional probability
$\Pcal[X_1(t)|X_1(0)]$. Thus, we can approximate the solution of
(\ref{th11-1}) by a stochastic process $Y_{t}$ along time intervals
$[t,t+\Delta t]$ using the following two-step algorithm for each
$\Delta t$, which we call {\bf splitting-step algorithm}:
\begin{enumerate}
\item[Step 1.] Knowing the value of $Y_t$ we obtain an intermediate value
$\tilde Y_t$ which is obtained through the exact integration of
(\ref{ss-1}) and $\tilde Y_t = X_1(t+\Delta t)$ and with initial
condition $X_1(t) = Y_t$.
\item[Step 2.] Then $\tilde Y_t$ is used as the initial condition for
(\ref{ss-2}) which is now integrated using any converging deterministic numerical
algorithm to get $\tilde X_2$. Then $Y_{t+\Delta t} = \tilde X_2(t+\Delta
t)$.
\end{enumerate}

The advantage of this splitting-step technique for SDEs subject to
boundary conditions is  is that if equation (\ref{ss-1}) is simple
enough and we know the solution $X_1$ of equation (\ref{th11-1}),
then the stochastic part of the problem can be handled correctly.
For example, for the case $\beta(X_1,t) = 0$ and $\sigma(X_1,t) =
\sqrt{X_1}$, it is known that the conditional probability
distribution is given by
\begin{eqnarray}
\Pcal[X_1(t)|X_1(0)]&=&\frac{2}{t}
\left(\frac{X_1(0)}{X_1(t)}\right)^{1/2} I_1 \left(\frac{4}{t}
\sqrt{X_1(t)X_1(0)}\right) e^{-\frac{2}{t} [X_1(t)+X_1(0)]} + \\
\nonumber & &  \qquad + e^{-\frac{4}{t} X_1(0)}\label{cdf-1}
\delta(X_1(t)),
\end{eqnarray}
where $I_1$ is the modified Bessel function and $\delta(x)$ is the
Dirac delta function.

This cdf can be sampled using the rejection or inverse methods but
this is computationally expensive. However, we introduce here a
very simple method \cite{moro} (see also \cite{broadie}) for
obtaining $X_1(t)$ by noting that the variable $Z(t) = \frac{4}{t}
X_1(t)$ has a probability distribution given by the non-central
$\chi^2$-distribution, that is
\begin{equation}\label{cdf-2}
\Pcal[Z|Z_0] = \sum_{j=1}^\infty \frac{(\lambda/2)^j
e^{-\lambda/2}}{j!} \Pcal_{\chi^2_{2 j}}(Z)+ e^{-\lambda/2}
\delta(Z)
\end{equation}
where $\lambda = \frac{4}{t} Z_0$ and $\Pcal_{\chi^2_{2j}}(X)$ is
the $\chi^2$-pdf with $2j$ degrees of freedom. Equation
(\ref{cdf-2}) reveals that the probability distribution for $Z$ is
a linear combination of $\chi^2$-pdfs with Poisson weights. This
fact can be exploited to generate $X_1(t)$ efficiently. If we
choose $K$ from a Poisson distribution with mean $\lambda/2$, then
\begin{equation}\label{rhodett}
X_1(t) = \frac{1}{2k} \left\{\begin{array}{ll}0 & \mathrm{if}\ K =
0,
\\ \sum_{i=1}^{2K} z_i^2 & \mathrm{if}\ K \neq 0, \end{array}\right.
\end{equation}
where $z_i$ are independent Gaussian random numbers with zero mean
and unit variance. Computationally, it is faster to sample the
random number $\sum_{i=1}^{2K} z_i^{2}$ using standard algorithms
for random number generation of the $\chi^2$ distribution. Other
examples of this sampling can be found in Section
\ref{sec-numerical} and the Appendix.

The obvious advantage of the splitting algorithm is that we
exploit the structure of the original equation more efficiently.
For example, in the above example, we get numbers from
(\ref{rhodett}) which are nonnegative and thus, if $\alpha(X,t)$
has nice properties, the approximated values for $X(t)$ are
nonnegative too. Such a splitting algorithm is not known from the
literature to the best of our current knowledge, although the idea
of splitting is not new for dynamical systems and their numerical
integration. A different kind of splitting technique has been
suggested in \cite{higham-mao-stuart}. However, their algorithm
called split-step Euler method is related to another subclass of
splitting of SDEs and their resulting split-step algorithm is of
lower order $0.5$ of mean square convergence, which is restricted
by their use of (drift-implicit) backward Euler method. Their
method indirectly refers to the splitting
\begin{eqnarray*}
d X_1 (t) & = & [ \alpha(X_1(t),t) + \beta(X_1(t),t) ] dt ,\\
d X_2 (t) & = & \sigma(X_2(t),t) dW(t)
\end{eqnarray*}
of the original system (\ref{inicio}) in our set-up,
where both equations for $X_1$ and $X_2$ are numerically integrated in a
separated fashion. In contrast to them, by a more efficient choice of splitting, we suggest
even to remove stochastic integrals from numerical integration by appropriate splitting
of (\ref{inicio}) in order to achieve higher order of convergence and our technique is not
only restricted to techniques of direct pathwise simulation. It may be noted that the explicit
removal of stochastic integrals from numerical integration by splitting
techniques is always possible and leads to converging numerical approximations
with higher order under some appropriate assumptions on the diffusion coefficient $\sigma$
(such as $\sigma \sigma \prime $ has sufficiently many bounded derivatives) in $\R^1$.
Another form of splitting has been suggested by \cite{castell-gaines}. They
present an algorithm which also takes advantage of techniques of numerical
integration of ODEs. However, both \cite{castell-gaines} and \cite{higham-mao-stuart}
do not discuss the issue of pathwise preservation of nonnegativity, monotonicity and
boundedness by their numerical approximation techniques.
%Habra que buscar si este tipo de algoritmos de splitting-step ya
%estan en la literatura.
%Ans: No estoy segurro en conjunto la literatura de algoritmos de splitting-step.

\section{General Convergence Theorems for Nonautonomous
Equations}\label{sec-teoremas} Let $C^{i,j}(\R^d \times [0,T])$
denote the vector space of continuous functions $f=f(x,t)$ which
are $i$ times continuously differentiable with respect to the
space-coordinate $x_k \in \R$ ($k=1,2,...,d$) and $j$ times
continuously differentiable with respect to time-coordinate $t \in
[0,T]$.

\subsection{A first general convergence theorem}
Recall that the original equation is
\begin{equation}\label{A1}
d X(t) = [\alpha(X(t),t) + \beta(X(t),t)] dt + \sigma(X(t),t) dW(t) .
\end{equation}
For the proof of splitting techniques below, we refer to the splitting
\begin{eqnarray}
\label{A2}
d X_1(t) & = & \beta(X_1(t),t) dt + \sigma(X_1(t),t) dW(t) , \\
\label{A3}
d X_2(t) & = & \alpha(X_2(t),t) dt  .
\end{eqnarray}

\begin{theorem}
\label{th:prop} Assume that the coefficient functions
$\alpha, \beta \in C^{2,1} (\R^d\times [0,T])$ and $\sigma \in C^{3,2}
(\R^d\times [0,T])$ with exclusively uniformly bounded derivatives are such that
$$\sup_{0 \le t \le T} \E \left[|X(t)|^2+|\alpha(X(t),t)|^2+|\beta(X(t),t)|^2
+|\sigma(X(t),t)|^2 \right] < + \infty$$ for a fixed finite, nonrandom
terminal time $T>0$.
Then the splitting-step algorithm with steps 1 and 2 has (global) strong and
weak order $1.0$ of convergence on the interval $[0,T]$ (in the worst case).
\end{theorem}

\begin{proof}
For simplicity, suppose that $d=1$. Let
$$0=t_0 < t_1 < ... < t_n < t_{n+1} < ... < t_N = T $$
be any nonrandom partition of the given time-interval $[0,T]$ with
sufficiently small maximum step size
$$\Delta \; = \; \max_{i=1,2,...,N} |t_i-t_{i-1}| \; \le \; 1. $$
Define the local pathwise error by
$\varepsilon^{loc}_{n+1} = X(t_n+\Delta_n t) - \tilde{X}(t_n+\Delta_n t) $
assuming that both exact solution $X$ and its approximation $\tilde{X}$ have
started at the same value $X(t_n)$ at time $t_n$.
To investigate this error, apply stochastic Taylor approximations to the
processes $X$, $X_1$ and $X_2$ -- an idea which originates from the Wagner-Platen
expansion \cite{wagner-platen78} and was popularized by \cite{sdenumeric0}
(more precisely speaking, this exploits the idea of an iterative
application of It\^{o} formula). For the sake of abbreviation,
we shall write $z_n$ or $z_u$ for all occurring coefficients or processes
(not referring to partial derivatives here), hence
$X_n=X(t_n)$, $\alpha_n=\alpha(X(t_n),t_n)$, $\beta_n=\beta(X(t_n),t_n)$,
$\sigma_n=\sigma(X(t_n),t_n)$, similarly $\alpha_u=\alpha(X(u),u)$,
$\beta_u=\beta(X(u),u)$, $\sigma_u=\sigma(X(u),u)$ and so on, unless it is
stated differently wherever convenient. Furthermore, define the partial
differential operators
\begin{eqnarray*}
L^0_0 f(x,t) & = & \frac{\partial f(x,t)}{\partial t} + [\alpha(x,t) + \beta(x,t)]
\frac{\partial f(x,t)}{\partial x} + \frac{1}{2} [\sigma(x,t)]^2
\frac{\partial^2 f(x,t)}{\partial x^2} ,\\
L^1_0 f(x,t) & = & \sigma(x,t) \frac{\partial f(x,t)}{\partial x} , \; L^1_1
f(x,t) \; = \; L^1_0 f(x,t) , \\
L^0_1 f(x,t) & = & \frac{\partial f(x,t)}{\partial t} + \beta(x,t)
\frac{\partial f(x,t)}{\partial x} + \frac{1}{2} [\sigma(x,t)]^2
\frac{\partial^2 f(x,t)}{\partial x^2} ,\\
L^0_2 f(x,t) & = & \frac{\partial f(x,t)}{\partial t} + \alpha(x,t)
\frac{\partial f(x,t)}{\partial x}
\end{eqnarray*}
where $L^0_i$ is mapping from $C^{2,1}(\R^d \times [0,T])$ to $C^{0,0}(\R^d \times [0,T])$
and $L^1_i$ from $C^{1,0}(\R^d \times [0,T])$ to $C^{0,0}(\R^d \times [0,T])$
for $i=1,2$, and $L^0_2$ from $C^{1,1}(\R^d \times [0,T])$ to $C^{0,0}(\R^d \times [0,T])$.
First, apply stochastic Taylor approximations to the solutions of (\ref{A1}) to
obtain
\begin{eqnarray}
\nonumber
X_{n+1} & = & X_n + \int^{t_{n+1}}_{t_n} [\alpha(X(s),s) + \beta(X(s),s)] ds +
\int^{t_{n+1}}_{t_n} \sigma(X(s),s) dW(s) \\
\nonumber
& = & X_n + [\alpha_n + \beta_n] \Delta_n t + \sigma_n \Delta_n W +\\
\nonumber
& & + \int^{t_{n+1}}_{t_n} \int^s_{t_n} L^0_0 (\alpha_u+\beta_u) \;
%\left[ \frac{\partial \alpha_u}{\partial u} + [\alpha_u+\beta_u] [\frac{\partial
%\alpha_u}{\partial x} + \frac{\partial \beta_u)}{\partial x} ] +
%\frac{1}{2} [\sigma_u]^2 [ \frac{\partial^2 \alpha_u}{\partial
%x^2}+ \frac{\partial^2 \beta_u}{\partial x^2} ] \right]
du \, ds + %\\ \nonumber & & +
\int^{t_{n+1}}_{t_n} \int^s_{t_n} L^1_0 (\alpha_u+\beta_u) \;
%\sigma_u \left[\frac{\partial \alpha_u}{\partial x} + \frac{\partial \beta_u}{\partial x} \right]
dW(u) \, ds + \\
\nonumber
& & + \int^{t_{n+1}}_{t_n} \int^s_{t_n} L^0_0 \sigma_u \;
%\left[ \frac{\partial \sigma_u}{\partial u} + [\alpha_u+\beta_u] \frac{\partial
%\sigma_u}{\partial x}  + \frac{1}{2} [\sigma_u]^2 \frac{\partial^2 \sigma_u}{\partial x^2} \right]
du \, dW(s) + %\\ & & +
\int^{t_{n+1}}_{t_n} \int^s_{t_n} L^1_0 \sigma_u \;
%\sigma_u \frac{\partial \sigma_u}{\partial x}
dW(u) \, dW(s) \\
& = & X_n + [\alpha_n + \beta_n] \Delta_n t + \sigma_n \Delta_n W +
\frac{1}{2} L^1_0 \sigma_n [(\Delta_n W)^2 - \Delta_n t] + R_{0,n}
\label{A4}
\end{eqnarray}
with remainder term
\begin{eqnarray*}
R_{0,n} & = & \int^{t_{n+1}}_{t_n} \int^s_{t_n} L^0_0 (\alpha_u+\beta_u) \; du \, ds +
\int^{t_{n+1}}_{t_n} \int^s_{t_n} L^1_0 (\alpha_u+\beta_u) \; dW(u) \, ds + \\
& & + \int^{t_{n+1}}_{t_n} \int^s_{t_n} L^0_0 \sigma_u \; du \, dW(s) +
\int^{t_{n+1}}_{t_n} \int^s_{t_n} \int^u_{t_n} L^0_0 L^1_0 \sigma_v \; dv \, dW(u) \, dW(s)
+ \\
& & + \int^{t_{n+1}}_{t_n} \int^s_{t_n} \int^u_{t_n} L^1_0 L^1_0 \sigma_v \; dW(v) \, dW(u) \, dW(s)
\end{eqnarray*}
Now, we have to compare (\ref{A4}) with what we get from the expansion of
the splitting method. In the latter case, by application of Wagner-Platen
expansion \cite{wagner-platen78} again, we arrive at
\begin{eqnarray}
\nonumber
X_{1,n+1} & = & X_n + \int^{t_{n+1}}_{t_n} \beta(X_1(s),s) \; ds +
\int^{t_{n+1}}_{t_n} \sigma(X_1(s),s) \; dW(s) \\
\nonumber
& = & X_n + \beta_n \Delta_n t + \sigma_n \Delta_n W +\\
\nonumber
& & + \int^{t_{n+1}}_{t_n} \int^s_{t_n} L^0_1 \beta_u \; du \, ds +
\int^{t_{n+1}}_{t_n} \int^s_{t_n} L^1_1 \beta_u \; dW(u) \, ds + \\
\nonumber
& & + \int^{t_{n+1}}_{t_n} \int^s_{t_n} L^0_1 \sigma_u \; du \, dW(s) +
\int^{t_{n+1}}_{t_n} \int^s_{t_n} L^1_1 \sigma_u \; dW(u) \, dW(s) \\
& = & X_n + \beta_n \Delta_n t + \sigma_n \Delta_n W + \frac{1}{2} L^1_1
\sigma_n [(\Delta_n W)^2 - \Delta_n t] + R_{1,n}
\label{AX1}
\end{eqnarray}
with the remainder term
\begin{eqnarray*}
R_{1,n} & = & + \int^{t_{n+1}}_{t_n} \int^s_{t_n} L^0_1 \beta_u \; du \, ds +
\int^{t_{n+1}}_{t_n} \int^s_{t_n} L^1_1 \beta_u \; dW(u) \, ds + \\
& & + \int^{t_{n+1}}_{t_n} \int^s_{t_n} L^0_1 \sigma_u \; du \, dW(s) +
\int^{t_{n+1}}_{t_n} \int^s_{t_n} \int^u_{t_n} L^0_1 L^1_1 \sigma_v \; dv \, dW(u) \, dW(s) +
\\
& & + \int^{t_{n+1}}_{t_n} \int^s_{t_n} \int^u_{t_n} L^1_1 L^1_1 \sigma_v \; dv \, dW(u) \, dW(s) .
\end{eqnarray*}
Here, the coefficients $\beta$ and $\sigma$ involved in the above integrals are evaluated at
the arguments $(X_1(u),u)$ and $(X_1(v),v)$, respectively.
Similarly, by deterministic Taylor expansion for the local approximation
of (\ref{A3}) in the framework of the splitting method, one gets to
\begin{eqnarray}
%\nonumber
\qquad X_{2,n+1} & = & X_{1,n+1} + \int^{t_{n+1}}_{t_n} \alpha(X_2(s),s) \; ds
\label{AX2} \; = \; X_{1,n+1} + \alpha (X_{1,n+1},t_n) \Delta_n t + R_{2,n}
\end{eqnarray}
with remainder term
\begin{eqnarray*}
R_{2,n} & = & \int^{t_{n+1}}_{t_n} \int^s_{t_n} L^0_2 \alpha_u \; du \, ds \\
& = & \int^{t_{n+1}}_{t_n} \int^s_{t_n} \left[\frac{\partial
\alpha(X_2(u),u)}{\partial u} + \alpha(X_2(u),u) \frac{\partial
\alpha(X_2(u),u)}{\partial x} \right] \; du \, ds .
\end{eqnarray*}
Here, the coefficients $\alpha$ involved in the above integrals are evaluated at
the arguments $(X_2(u),u)$. An expansion of $\alpha(X_{1,n+1},t_n)$ with respect
to space-variable $x$ gives
\begin{eqnarray}
\qquad \alpha(X_{1,n+1},t_n) & = & \alpha_n + \int^{t_{n+1}}_{t_n} L^0_1
\alpha(X_1(s),s) \; ds + \int^{t_{n+1}}_{t_n} L^1_1 \alpha(X_1(s),s) \; dW(s) .
\label{Aalpha}
\end{eqnarray}
Now, plug expansions (\ref{Aalpha}) and (\ref{AX1}) into the expansion
(\ref{AX2}) in order to encounter with
\begin{eqnarray}
\label{A5}
X_{2,n+1} & = & X_n + [\alpha_n + \beta_n] \Delta_n t + \sigma_n \Delta_n W + \frac{1}{2} L^1_1
\sigma_n [(\Delta_n W)^2 - \Delta_n t] + R_{3,n}
\end{eqnarray}
with the remainder term
\begin{eqnarray*}
R_{3,n} & = & R_{1,n} + R_{2,n} + \left[\int^{t_{n+1}}_{t_n} L^0_1 \alpha_s
ds + \int^{t_{n+1}}_{t_n} L^1_1 \alpha_s dW(s) \right] \Delta_n t  .
\end{eqnarray*}
Consequently, the local pathwise error $\varepsilon^{loc}_{n+1} = X_{n+1} -
X_{2,n+1}$ can be represented by
\begin{eqnarray}
\nonumber
\varepsilon^{loc}_{n+1} & = & X_{n+1} - X_{2,n+1} \; = \; R_{0,n} - R_{3,n} \\
%& = & [\alpha' \beta + \frac{1}{2}\alpha''\sigma^2]
%\Delta_n t + \alpha' \sigma \Delta_n t \Delta_n W +\alpha' \sigma
%[\alpha''\sigma + \alpha'\sigma'] \Delta_n t [\Delta_n W - \Delta_n t] \nonumber \\
%& & - \alpha'\sigma \Delta_n Z - \frac{1}{2}[\beta\alpha'+\alpha\beta'+\frac{1}{2}\sigma
%\alpha'']\Delta_n^2 - \alpha \sigma' [\Delta_n W \Delta_n t - \Delta_n Z]
\nonumber
& = & \int^{t_{n+1}}_{t_n} \int^s_{t_n} L^0_0 (\alpha_u+\beta_u) \; du \, ds +
\int^{t_{n+1}}_{t_n} \int^s_{t_n} L^1_0 (\alpha_u+\beta_u) \; dW(u) \, ds + \\
\nonumber
& & + \int^{t_{n+1}}_{t_n} \int^s_{t_n} L^0_0 \sigma_u \; du \, dW(s) +
\int^{t_{n+1}}_{t_n} \int^s_{t_n} \int^u_{t_n} L^0_0 L^1_0 \sigma_v \; dv \, dW(u) \, dW(s)
+ \\
\nonumber
& & + \int^{t_{n+1}}_{t_n} \int^s_{t_n} \int^u_{t_n} L^1_0 L^1_0 \sigma_v \; dW(v) \, dW(u) \, dW(s) +
\\
\nonumber
& & - \int^{t_{n+1}}_{t_n} \int^s_{t_n} L^0_1 \beta_u \; du \, ds -
\int^{t_{n+1}}_{t_n} \int^s_{t_n} L^1_1 \beta_u \; dW(u) \, ds + \\
\nonumber
& & - \int^{t_{n+1}}_{t_n} \int^s_{t_n} L^0_1 \sigma_u \; du \, dW(s) -
\int^{t_{n+1}}_{t_n} \int^s_{t_n} \int^u_{t_n} L^0_1 L^1_1 \sigma_v \; dv \, dW(u) \, dW(s) +
\\
\nonumber
& & - \int^{t_{n+1}}_{t_n} \int^s_{t_n} \int^u_{t_n} L^1_1 L^1_1 \sigma_v \; dv \, dW(u) \, dW(s)
- \int^{t_{n+1}}_{t_n} \int^s_{t_n} L^0_2 \alpha_u \, du \, ds +
\\
& & - \left[\int^{t_{n+1}}_{t_n} L^0_1 \alpha_s \;
ds + \int^{t_{n+1}}_{t_n} L^1_1 \alpha_s \; dW(s) \right] \Delta_n t  .
\label{A6}
\end{eqnarray}
%where $(\Delta_n Z,\Delta_n W)$ are Gaussian distributed and the coefficients $\alpha, \beta$ and
%$\sigma$ are calculated at (random) intermediate points.
Now recall that $\alpha, \beta, \sigma$ have exclusively uniformly bounded
derivatives and their second moments along the solution of (\ref{A1}) are
bounded on $[0,T]$. Therefore, all operators $L^j_i$ applied to coefficients
$\alpha, \beta$ and $\sigma$ have images with uniformly bounded second
moments. This implies that $\E\varepsilon^{loc}_n =
\Ocal([\Delta_n t]^2)$ and $\E (\varepsilon^{loc}_n- \E\varepsilon^{loc}_n)^2 = \Ocal([\Delta_n t]^{3})$.
For mean square convergence, it remains to apply fairly general convergence theorems as known from
\cite{sdenumeric1} or \cite{sdenumeric4} with local rates $r_1 = 2$ and $r_2=1.5$ in order to conclude
the global strong (i.e. in the $L^2$-sense) convergence order $r=1.0$. More
precisely, we have
$$ \max_{n = 0, 1, ..., N} \left(\E |X(t_n) - \tilde{X}(t_n)|^2 \right)^{1/2} \; =
\; \Ocal (\Delta) $$
provided that the initial errors $\E |X(0) - \tilde{X}(0)|^2 = \Ocal
(\Delta^2) $ started at $L^2$-integrable initial values which are independent of the
$\sigma$-algebra generated by the underlying Wiener process $W$.
Eventually, the global order $r_w=1.0$ of weak convergence with respect to smooth,
polynomially bounded test functions is rather obvious from the fact that the local weak rate
$r_w=2.0$ (exploiting standard techniques known from \cite{feller} and \cite{sdenumeric1}
on weak convergence analysis). Consequently, the splitting-step algorithm has the claimed
convergence orders under the above stated assumptions, and the proof is complete.
\qquad\end{proof}

\subsection{Numerical generation of strong path solutions for some
processes}
The success of the splitting-step method requires the
exact numerical integration of step 1 which can be problematic in
some specific situations. That is why we have applied in previous
numerical examples such equations which allow us to represent the
solution from step 1 in terms of pure functions of the Wiener
process and initial values such as $X(t)=F(X(0),W(t))$ or simple
deterministic integrals as seen with the example of Bessel-type
diffusions $X(t)=(\sqrt{X(0)}+W(t))^2$ or the geometric Brownian
motion $X(t)=X(0)\exp((\alpha-\sigma^2/2) t + \sigma W(t))$. Also,
for the Ornstein-Uhlenbeck processes or logistic equations, one
may use well-known integration-by-parts formula and / or the
information on the exact distribution of involved stochastic
integrals with deterministic differentials which can be pathwisely
treated by deterministic quadrature methods. As another
alternative, one could exploit the Doss representation (see
\cite{doss}) to develop a semi-analytic approach or a ODE-PDE
approach to decompose the original problem into a random ODE and a
deterministic PDE in order to figure out which numerical
implementation is more efficient in conjunction with our splitting
technique. However, in general one produces additional
discretization errors which might influence the accuracy of the
computations using the splitting algorithm as well. Then it is a
must to consider its stability properties too. Such a fairly
complex and very problem-dependent work we leave to the future.

\subsection{A general theorem on $L^2$-convergence based on VOP}
To relax some of the technical assumptions, once may also exploit the
variation-of-constants formula (VOP). Suppose that the original equation is
\begin{equation}\label{S0}
\quad \; \; d X(t) = [\alpha(X(t),t) + \beta(X(t),t)] dt + \sigma_1(X(t),t) dW_1(t) +
\sigma_2(X(t),t)) dW_2(t)
\end{equation}
where $W_1$ and $W_2$ are independent Wiener processes.
Consider the splitting
\begin{eqnarray}
\label{S1}
d X_1(t) & = & \beta(X_1(t),t) dt + \sigma_1(X_1(t),t) dW_1(t) , \\
\label{S2}
d X_2(t) & = & \alpha(X_2(t),t) dt  + \sigma_2(X_2(t),t) dW_2(t) .
\end{eqnarray}
Let $C^{0}_{locLip}(S)$ denote the vector space of local Lipschitz
continuous functions on the open set $S$.

\begin{theorem}
\label{th:sprop} Assume that the coefficient functions
$\alpha, \beta \in C^{0}_{locLip} (\D\times [0,T])$ and $\sigma_i \in
C^{0}_{locLip} (\D\times [0,T])$ are such that the continuous unique strong solution $X$
of (\ref{S0}) exists on the closed set $\bar{\D}$,
$$ \sup_{0 \le t \le T} \E \Big[|X(t)|^2\Big] +
\sup_{0 \le t < T} \sup_{t \le s \le T} \E \Big[|\Phi(s,X(t))|^2\Big] < + \infty$$
$$
\sup_{0 \le t \le T} \E \left[ \int^T_t \!\!\!|[\Phi(s,X(t))]^{-1}\alpha(X(s),s)|^2 ds
+ \!\int^T_t \!\!\!|[\Phi(s,X(t))]^{-1}\sigma_2(X(s),s)|^2 ds \right]
%|\alpha(X(t),t)|^p+|\beta(X(t),t)|^p +|\sigma_1(X(t),t)|^p + |\sigma_2(X(t),t)|^p
\!<\! + \infty$$
for a fixed finite, nonrandom terminal time $T>0$ %%and sufficiently large exponents $p \ge 2$
and the stochastic flow $\Phi$ generated by (\ref{S1}) is mean square
H\"older-continuous with exponent $r_H \ge 0.5$. Furthermore, suppose that
step 1 of the splitting algorithm is exactly integrable and step 2 can be
carried out with local mean accuracy with rate $r_1 \ge 1.0$, local
mean square accuracy with rate $r_2\ge 0.5$ and
$$min(r_1,r_H+1.0) \ge min(r_2,r_H+1.0) + 0.5 .$$
Then, in the worst case, the error of splitting-step algorithm with steps 1 and 2 has
(global) %strong and weak
order
$$r_g \ge min (1.0,min(r_2,r_H+1.0)-0.5)$$
of $L^2$-convergence on the interval $[0,T]$.
\end{theorem}

\begin{proof}
Suppose that (\ref{S1}) has known fundamental solution $\Phi=\Phi(t,X_0)$
with a.s. H\"older exponent $r_H$. Then the exact solution
of the original equation (\ref{S0}) possesses the pathwise
representation
\begin{eqnarray*}
X(t+\Delta t) & = &
\Phi(t+\Delta t,X(t)) + \Phi(t+\Delta t,X(t)) \int^{t+\Delta t}_t
[\Phi(s,X(t))]^{-1} \alpha (X(s),s) \, ds \\
& & + \Phi(t+\Delta t,X(t)) \int^{t+\Delta t}_t
[\Phi(s,X(t))]^{-1} \sigma_2 (X(s),s) \, dW_2(s)
%& \approx & \Phi(t+\Delta t,X(t)) + \alpha(X(t),t) \Delta t
\end{eqnarray*}
on each subintervals $[t,t + \Delta t] \subset [0,T]$.
Note that $\Phi(t+\Delta t,X(t))$ is independent of $W_2(s)-W_2(t)$ for $s
\ge t$. Let $Y(t)$ denote the value of right-continuous numerical approximation
of splitting step algorithm (which we always can construct using step functions
in a standard way). The main idea is to apply the fairly general $L^2$-convergence
theory known from \cite{sdenumeric1}, \cite{sdenumeric3} and \cite{sdenumeric4}.
For this purpose, we need to study the local conditional
accuracy of our splitting algorithm. Locally, we may suppose that
$X(t)=Y(t)=x$ is deterministic (${\mathcal F}_t$-adapted). First, consider the
local conditional mean accuracy of the splitting algorithm. We find that
\begin{eqnarray*}
\lefteqn{\left|\E \Big[ X(t+\Delta t) - Y(t+\Delta t)\Big|{\mathcal F}_t\Big]\right|}\\
& = & \left| \E \Big[ \Phi(t+\Delta t,x) \int^{t+\Delta t}_t
[\Phi(s,x)]^{-1} \alpha (X(s),s) \, ds + \right. \\
& & \left. + \Phi(t+\Delta t,x) \int^{t+\Delta t}_t [\Phi(s,x)]^{-1} \sigma_2 (X(s),s) \, dW_2(s)
- (Y(t+\Delta t)-x) \Big| {\mathcal F}_t \Big] \right|\\
& \le &
\left| \E \Big[ \int^{t+\Delta t}_t \Big(\Phi(t+\Delta t,x)-\Phi(s,x)\Big) [\Phi(s,x)]^{-1}
\alpha (X(s),s) \, ds \Big| {\mathcal F}_t \Big] \right| + \\
& & + \left| \E \Big[ \int^{t+\Delta t}_t \Big(\Phi(t+\Delta t,x)-\Phi(s,x)\Big) [\Phi(s,x)]^{-1}
\sigma_2 (X(s),s) \, dW_2(s) \Big| {\mathcal F}_t \Big] \right| + \\
& & + \left| \E \Big[\int^{t+\Delta t}_t \alpha (X(s),s) ds + \int^{t +
\Delta t}_t \sigma_2 (X(s),s) dW_2(s) - (Y(t+\Delta t)-x) \Big| {\mathcal F}_t \Big]
\right|\\
& \le & K_1 [\Delta t]^{min(r_1,(2r_H+1)/2+0.5)} ,
\end{eqnarray*}
where $K_1$ is a real constant, hence the local conditional mean accuracy rate
$$min(r_1,(2r_H+1)/2+0.5) = min(r_1,r_H+1.0)$$
of the splitting algorithm can be verified. Second, consider the local conditional mean
square accuracy of the splitting algorithm.
\begin{eqnarray*}
\lefteqn{\left(\E \Big[ |X(t+\Delta t) - Y(t+\Delta t)|^2 \Big|{\mathcal
F}_t\Big]\right)^{1/2}}\\
& = & \left( \E \Big[ \Big| \Phi(t+\Delta t,x) \int^{t+\Delta t}_t
[\Phi(s,x)]^{-1} \alpha (X(s),s) \, ds + \right. \\
& & \left. + \Phi(t+\Delta t,x) \int^{t+\Delta t}_t [\Phi(s,x)]^{-1} \sigma_2 (X(s),s) \, dW_2(s)
- (Y(t+\Delta t)-x) \Big|^2 \Big| {\mathcal F}_t \Big] \right)^{1/2}\\
& \le &
\left( \E \Big[ \Big|\int^{t+\Delta t}_t \Big(\Phi(t+\Delta t,x)-\Phi(s,x)\Big) [\Phi(s,x)]^{-1}
\alpha (X(s),s) \, ds \Big|^2 \Big| {\mathcal F}_t \Big] \right)^{1/2} + \\
& + & \left( \E \Big[ \Big| \int^{t+\Delta t}_t \Big(\Phi(t+\Delta t,x)-\Phi(s,x)\Big) [\Phi(s,x)]^{-1}
\sigma_2 (X(s),s) \, dW_2(s) \Big|^2 \Big| {\mathcal F}_t \Big]
\right)^{1/2} \!\!+ \\
& + & \left( \E \Big[ \Big| \int^{t+\Delta t}_t \!\!\!\!\!\alpha (X(s),s) ds + \int^{t +
\Delta t}_t \!\!\!\!\!\sigma_2 (X(s),s) dW_2(s) - (Y(t+\Delta t)-x) \Big|^2 \Big| {\mathcal F}_t \Big]
\right)^{1/2}\\
& \le & K_2 [\Delta t]^{min(r_2,(2r_H+1)/2+0.5)} ,
\end{eqnarray*}
where $K_2$ is a real constant, hence the local conditional mean square accuracy rate
$$min(r_2,(2r_H+1)/2+0.5) = min(r_2,r_H+1.0)$$
of the splitting algorithm can be derived. Now, apply the general
convergence theory from \cite{sdenumeric1}, \cite{sdenumeric3} and \cite{sdenumeric4} to
conclude the worst case estimate of global mean square convergence rate $r_g$ as
stated in Theorem \ref{th:sprop} along nonrandom partitions of $[0,T]$. This
completes the proof of Theorem \ref{th:sprop}.
\end{proof}

\vspace*{5mm}

\noindent
{\bf Remark.}
Consider the example
\begin{equation}
dX(t) = (1+X(t)) dt + 2 \sqrt{X(t)} dW(t)
\end{equation}
with $\alpha(X,t) = X$, $\beta(X,t) = 1$ and $\sigma_1(X,t) =
2\sqrt{X}$ and $\sigma_2(X,t)=0$.  This example satisfies the
assumptions of Theorem \ref{th:sprop} with $r_H=0.5$. To see this fact, one
has to show that the related stochastic flow is given by
the Bessel-type flow
$$ \Phi(t,X(0)) = (\sqrt{X(0)}+W(t))^2 $$
while using It\^o formula. Moreover, this flow has uniformly bounded
moments of all orders for nonrandom initial data and
%and an inverse $$ [\Phi(s,Z(t))]^{-1} = [\sqrt{Z(t)}-(W(s)-W(t))]^2 $$
%with uniformly bounded moments too.
is mean square H\"older-continuous with exponent $r_H=0.5$ since
$$ \Phi(t,X(0))-\Phi(s,X(0)) = 2 \sqrt{X(0)} (W(t)-W(s)) $$
and
$$ \E |\Phi(t,X(0))-\Phi(s,X(0))|^2 = 2 \E [X(0)] \E |W(t)-W(s)|^2 = 2 \E
[X(0)] |t-s| $$
provided that $X(0) \ge 0$ (a.s.) is independent of the process $W$. Therefore, the
related splitting-step algorithm can achieve a global $L^2$-convergence order
$r_g = 1.0 $ since the numerical integration of step 2 can be implemented by
any deterministic Runge-Kutta method with an interplay of local accuracy
rates $r_1 = 2.0$ and $r_2 = 1.0$. Similarly, we can proceed for other
equations with $\sqrt{(.)}$- or other H\"older-continuous terms with
H\"older exponent $\ge 0.5$.

\section{First illustrative numerical experiments}
\label{sec-numerical} In this section we will give several
illustrative and important examples of applying our new stochastic
numerical schemes to SDEs and test the rate of convergence
obtained in previous Section \ref{sec-teoremas}.
\subsection{Using the transition probability}
\begin{figure}
\begin{center}
\includegraphics[width=3.0in,clip=]{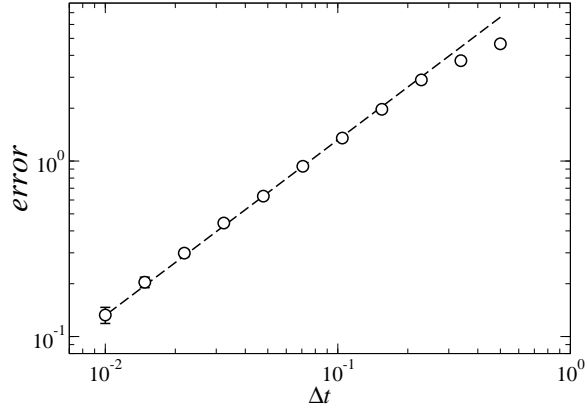}
\caption{\label{fig-error} Absolute value of the error
(\ref{eps1}) as a function of $\Delta t$ for $t = 1.0$ and $X(0) =
1$ obtained using the splitting-step algorithm. The dashed line is
proportional to $\Delta t$.}
\end{center}
\end{figure}

We applied the splitting-step method to the following SDE
\begin{equation}
dX(t) = (1+X(t)) dt + 2 \sqrt{X(t)} dW(t)
\end{equation}
with $\alpha(X,t) = 1 + X$, $\beta(X,t) = 0$ and $\sigma(X,t) =
2\sqrt{X}$. The conditional mean value is given by
\begin{equation}
\E (X(t)|X(0)=x_0) = (x_0 + 1) e^t - 1
\end{equation}
for any non-random value $x_0 \ge 0$. In figure \ref{fig-error}
the error
\begin{equation}\label{eps1}
\varepsilon_1 = | \E \tilde X(t) - [(x_0)e^t -1]|
\end{equation}
versus decreasing uniform step size $\Delta t$ is depicted,
where $\tilde X(t)$ is the solution obtained through the numerical
approximation.  Figure \ref{fig-error} shows statistical-numerical
evidence that the method has weak order 1.0 while using constant
step sizes $\Delta t$.

\subsection{Using the exact solution}
\begin{figure}
\begin{center}
\includegraphics[width=3.0in,clip=]{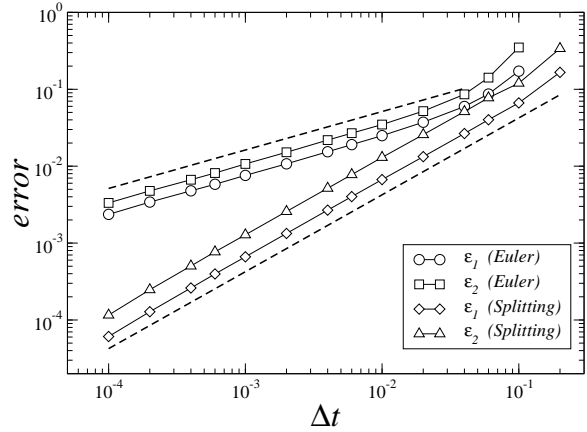}
\caption{\label{fig-error1} Value of the error (\ref{eps2}) as a
function of $\Delta t$ for $t = 5.0$ and $X(0) = 1.0$ using the
splitting-step and Euler algorithms. Dashed lines are the power
laws $\Delta t$ and $\sqrt{\Delta t}$.}
\end{center}
\end{figure}

We now apply the splitting-step algorithm to the stochastic
Ginzburg-Landau equation
\begin{equation}\label{exp2}
dX(t) = (X(t)-[X(t)]^3) \, dt + X(t) \, dW(t) .
\end{equation}
In this case we take advantage of the known exact solution for the
linear part of this equation which is
\begin{equation}\label{exp2-2}
\qquad dX_1(t) = X_1(t) dt + X_1(t) dW(t) \Rightarrow X_1(t) = X_1(s) \exp \Big(
(t-s)/2 + W(t) - W(s)\Big) \!.
\end{equation}
To integrate the remaining part of the equation we use an
partial-implicit nonstandard technique which is a
nonnegativity-preserving one\footnote{This is true only if
$[X_2(0)] \sup_{n\in N} \Delta_n t < 2$.} given by
\begin{equation}\label{exp2-3}
X_2(t+\Delta t) = X_2(t) - \frac{\Delta t}{2} X_2^2 (t) [X_2(t)
+X_2(t+\Delta t)]
\end{equation}
Our results for the strong error
\begin{equation}\label{eps2}
\varepsilon_k(t) = ( \E| X(t)-\tilde X(t)|^k)^{1/k}
\end{equation}
are shown in Figure \ref{fig-error1} for $k=1,2$ and compared with
the Euler algorithm for the same equation. Obviously, our results indicate
that our proposed splitting method indeed is of strong order 1.0 while using
constant step sizes $\Delta t$.

\section{Stochastic models in finance}
In this section we apply the splitting-step method to some
fundamental models in mathematical finance. This application will
also serve to introduce other algorithmically simple samplings of
conditional probability transitions as for the $\sqrt{X_1(t)}$
case shown in the introduction.

\subsection{Interest rate model of Cox-Ingersoll-Ross}
\begin{figure}
\begin{center}
\includegraphics[width=3.0in,clip=]{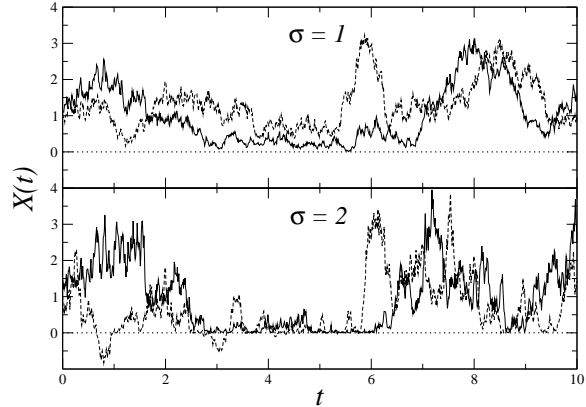}
\caption{\label{figcir} Numerical approximation of the CIR process
(\ref{coxonceagain}) using the splitting-step method
(\ref{ssCIR}), (solid line) and the method proposed in
\cite{higham-mao} (dashed line) with $a=b=1$ and $\Delta t=
10^{-2}$ and $X(0) = 1$ and different values of $\sigma$: Upper
panel shows the case $a \geq \sigma^2/2$ where the boundary at
zero becomes unattainable, while lower panel is for $a <
\sigma^2/2$. The splitting-step numerical approximation reproduces
the attainability of the boundary condition while preserving the
non-negativity of the solution of (\ref{coxonceagain}), whereas
the method in \cite{higham-mao} has trajectories which take
negative values for $X(t)$.}
\end{center}
\end{figure}

An interesting example in which the splitting-step scheme is
particularly efficient is the Cox-Ingersoll-Ross (CIR) model
\cite{cox} for stochastic interest rates
\begin{equation}\label{coxonceagain}
dX(t) = [a+b X(t)] dt + \sigma \sqrt{X(t)} dW(t), \quad X(0) = x_0
\geq 0
\end{equation}
with real parameters $a\geq 0$, $b \in \R$ and $\sigma > 0$. As
stated in the introduction, strong solutions of
(\ref{coxonceagain}) are nonnegative. However, depending on the
specific values of parameters $a$, $b$ and $\sigma$,
distinct behavior at the boundary is possible.
\begin{itemize}
\item If $a\geq \sigma^2/2$ then the solution is always positive
$X(t) > 0$ if $x_0>0$, because the boundary $X_(t)=0$ becomes
unattainable. \item If $a < \sigma^2/2$ there are infinite many
values of $t>0$ for which $X(t) =0$. The boundary becomes
attainable, but it is (instantaneously) reflecting. That is, when
a sample path reaches 0, then it returns immediately to the
interior of the state space in a reflecting manner.
\end{itemize}
The exact transition density for the CIR process is known, but its
sampling can be difficult depending on the parameters $a$, $b$ and
$\sigma$. Here we may exploit the simplest splitting which appropriately
reflects the boundary behavior of the CIR process by SDEs
\begin{eqnarray}\label{ssCIR}
d X_1(t) & = & a \, dt + \sigma\sqrt{X_1(t)} dW(t) \\
\label{ssCIR-1} d X_2(t)&=&b X_2(t) dt
\end{eqnarray}
which can be easily inferred by noting that the boundary behavior
does not depend on the parameter $b$. To integrate the first step
in the splitting-step system (\ref{ssCIR})-(\ref{ssCIR-1}) we note
that the process defined by (\ref{ssCIR}) represents an
$a$-dimensional Bessel process \cite{feller}. Its transition
density $\Pcal[X_1(t+\Delta t)|X_1(t)]$ can be written in terms of
a non-central $\chi^2$ distribution and in particular, we have
that:
\begin{equation}\label{ssCIR1}
X_1(t+\Delta t) = \frac{\sigma^2\Delta t}{4} \chi'^2_d(\lambda)
\end{equation}
where
\begin{equation}\label{ssCIR2}
\lambda = \frac{4 X_1(t)}{\sigma^2 \Delta t},\qquad d =
\frac{4a}{\sigma^2}
\end{equation}
and $\chi'^2_d(\lambda)$ random numbers can be sampled using the
algorithms in the appendix. The last part of the splitting-step
scheme (\ref{ssCIR-1}) can be integrated exactly or using
numerical approximations. In our case we use deterministic Euler
approximations where non-negativity is preserved if $\Delta t$ is
small enough. Our simulations for the CIR process are shown in
Figure \ref{figcir}. The boundary behavior for different values of
$a$ and $\sigma$ is reproduced and, at the same time,
non-negativity is conserved. Note that other usual integration
strategies as that of \cite{higham-mao} eventually produce
negative values for $X(t)$ (for finite $\Delta t$) which lack any
possible interpretation in the context of finance and could induce
severe errors in option valuation.

\subsection{Constant Elasticity Volatility models}
Another important stochastic process in finance is the constant
elasticity of variance (CEV) diffusion to model asset prices. This
process, first introduced to finance by Cox \cite{cox1}, is
capable of reproducing the volatility smile observed in the
empirical data unlike other standard price models like the
Black-Scholes-Merton geometric Brownian motion. The process is
defined as
\begin{equation}\label{CEV}
dX(t) = \mu X(t) dt + \sigma X(t)^\gamma dW(t),\quad X(0) = x_0.
\end{equation}
The CEV model includes the geometric Brownian model of Black,
Scholes and Merton ($\gamma = 1$) and the square-root models of
Cox and Ross ($\gamma = 1/2$). Contrary to the Black, Scholes and
Merton model, the CEV model incorporates a variance adjustment
that causes the absolute level of the variance to decline as the
stock price rises and to rise as the stock price declines, as seen
empirically in most equity and interest rate volatilities.

For our purposes, the CEV model encompasses most of the boundary
conditions we can implement at $X(0)=0$. Depending on the value of
$\gamma$ the boundary $X(0) = 0$ is
\begin{itemize}
\item {\em Natural boundary} for $\gamma \geq 1$, which means that
boundary is unattainable in finite time. \item {\em Exit or
absorbing boundary} for $1/2\leq\gamma <1$, i.e. $X(t)$ reaches
zero with finite probability in finite time and gets absorbed at
it. \item {\em Regular boundary} for $\gamma < 1/2$. Now the
boundary can be reached in finite time, and we need to specify a
boundary condition at $X(0)=0$. Typical choices are reflecting or
absorbing ones.
\end{itemize}
The transition density for the CEV process is known for general
values of $\gamma$ and $\mu$ (see \cite{cox1}). However, due to
the fact that the boundary behavior of $X(t)$ does not depend on
$\mu$, we propose here the following splitting-step system
\begin{eqnarray}\label{ssCEV1}
dX_1&=&\sigma X_1^\gamma dW(t), \\
\label{ssCEV2} dX_2&=&\mu X_2 dt.
\end{eqnarray}
The classification of the boundary $X_1(t) = 0$ for (\ref{ssCEV1})
is the same as for (\ref{CEV}). We discuss now the sampling of the
transition density for each value of $\gamma$.\footnote{We do not
consider the case $\gamma =1$, since this is trivially integrated
using the strong solution (\ref{exp2-2}).}

\subsubsection{$\gamma > 1$, Natural Boundary} In this case
$X_1(t)=0$ is unattainable and its probability transition density
is given by
\begin{equation}\label{pdfssCEV1ag1}
\Pcal[x_t|x_0]=\frac{x_t^{1/2-2\gamma}x_0^{1/2}}{\sigma^2 (\gamma
-1)} \exp
\left[-\frac{x_0^{2(1-\gamma)}+x_t^{2(1-\gamma)}}{2\sigma^2(1-\gamma)^2
t}\right]
I_{\frac{1}{2(\gamma-1)}}\left[\frac{x_0^{1-\gamma}x_t^{1-\gamma}}{\sigma^2
(1-\gamma)^2 t}\right]
\end{equation}
To sample this distribution, we can consider the relationship
between the CEV and the CIR processes. Since $X_1(t)=0$ is not
accessible, we can make the following change of variables $Y_1 =
X_1^{2(1-\gamma)}/[4(\gamma -1)^2\sigma^2]$ to find that $dY_1 =
\lambda dt+ \sqrt{Y_1} dW(t)$ where $\lambda =
(1-2\gamma)/[4(1-\gamma)]$. Thus $Y_1$ is a Bessel process which
can be sampled using the non-central $\chi^2$ distribution. In the
end we have that
\begin{equation}\label{samplingCEVagt1}
X_1(t+\Delta t) = \left[(\gamma -1)^2\sigma^2 \Delta t\,
\chi'^2_d\left(\frac{X_1(t)^{2(1-\gamma)}}{\sigma^2(\gamma-1)^2\Delta
t}\right)\right]^{\frac{1}{2(1-\gamma)}}
\end{equation}
where $d=(1-2\gamma)/(1-\gamma) \geq 1/2$. A typical path of the
CEV process for $\gamma > 1$ is shown in Figure \ref{figcev}.

\subsubsection{$1/2\leq\gamma < 1$, Absorbing boundary} In this
case $X_1(t)=0$ is an absorbing boundary. The transition
probability density was found by Cox \cite{cox1} and is given by
\begin{equation}\label{pdfssCEV1ag12}
\Pcal[x_t|x_0]=\frac{x_t^{1/2-2\gamma}x_0^{1/2}}{\sigma^2
(1-\gamma)} \exp
\left[-\frac{x_0^{2(1-\gamma)}+x_t^{2(1-\gamma)}}{2\sigma^2(1-\gamma)^2
t}\right]
I_{\frac{1}{2(1-\gamma)}}\left[\frac{x_0^{1-\gamma}x_t^{1-\gamma}}{\sigma^2
(1-\gamma)^2 t}\right].
\end{equation}
This transition probability does not integrate to one because
there is a finite probability that the trajectory gets absorbed at
zero given by \cite{cox1}
\begin{equation}\label{absprobability}
\Pcal[0|X_1(t)] =
G\left(\frac{1}{2(1-\gamma)},\frac{x_1(t)^{-2(1-\gamma)}}
{2\sigma^2(1-\gamma)^2\Delta t}\right)
\end{equation}
where $G(\nu,x)$ is the complementary Gamma function. For this
equation, we know that trajectories of (\ref{ssCEV1}) get absorbed
at zero with probability one while taking the limit $\Delta t \to
\infty$. Actually, this is also the case for the solutions of
(\ref{CEV}) for $1/2\leq\gamma < 1$. This is counterintuitive
since the mean $\E [X_1(t)]$ is constant for (\ref{ssCEV1}) or
grows exponentially like $x_0 e^{\mu t}$ for (\ref{CEV}).  This
intriguing feature of (\ref{CEV}) hampers the numerical
simulations of this process and in fact only exact sampling of the
probability transition density (\ref{pdfssCEV1ag12}) and
(\ref{absprobability}) correctly accounts for it at finite $\Delta
t$.

\begin{figure}
\begin{center}
\includegraphics[width=3.0in,clip=]{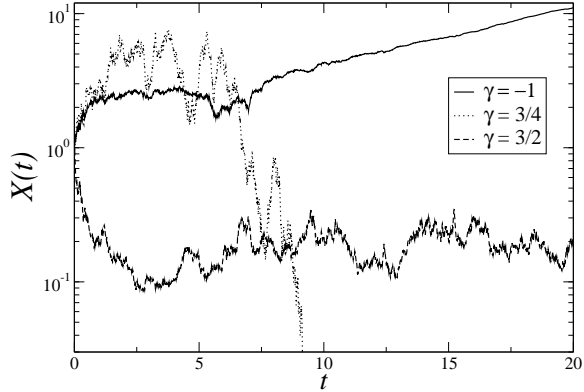}
\caption{\label{figcev} Numerical approximation of the CEV process
(\ref{CEV}) using the splitting-step method
(\ref{ssCEV1})-(\ref{ssCEV2}), for different values of $\gamma$.
The case $\gamma = -1$ is the solution of (\ref{CEV}) with
reflecting boundary conditions, while the case $\gamma = 3/4$ gets
absorbed at zero at finite time (not shown in the logarithmic
scale). Parameters are $\mu = 0.1$, $\Delta t = 10^{-2}$.}
\end{center}
\end{figure}

While similar to (\ref{pdfssCEV1ag1}), the transition probability
density (\ref{pdfssCEV1ag12}) cannot be sampled using the
relationship to the CIR process, since the solution of
$dX_1=X_1^\gamma dW(t)$ has a finite probability to be absorbed at
$X_1(t) = 0$. However, using the relationship $I_n(x) =
I_{-n}(x)$, $n\in \N$ for the modified Bessel function we have
that if $1/[2(1-\gamma)]=n$, $X_1$ can be sampled from a
non-central $\chi^2$ distribution with negative (integer) number
of degrees of freedom (see appendix):
\begin{equation}\label{samplingCEVagt12}
X_1(t+\Delta t) = \left[(\gamma -1)^2\sigma^2 \Delta t\,
\chi'^2_d\left(\frac{X_1(t)^{2(1-\gamma)}}{\sigma^2(\gamma-1)^2\Delta
t}\right)\right]^{\frac{1}{2(1-\gamma)}}
\end{equation}
where $d=2-2n=0,-2,-4,\ldots...$ and then $\gamma = 1-1/2n =
1/2,3/4,5/6,\ldots$. In figure \ref{figcev} we show a typical path
for the CEV process for $\gamma = 3/4$ ($d=-2$) in which we see
that it gets absorbed at zero for finite time.

\subsubsection{\bf $ 0 \leq \gamma < 1/2$, Regular Boundary}
Now the boundary is attainable and, contrary to the previous case,
a boundary condition should be provided. If we consider the
simplest cases, namely absorbing or reflecting boundary condition
then the transition probability density is given by
\begin{equation}\label{pdfssCEV1al12}
\Pcal_{\pm}[x_t|x_0]=\frac{x_t^{1/2-2\gamma}x_0^{1/2}}{\sigma^2
|1-\gamma|} \exp
\left[-\frac{x_0^{2(1-\gamma)}+x_t^{2(1-\gamma)}}{2\sigma^2(1-\gamma)^2
t}\right]
I_{\pm\frac{1}{2|1-\gamma|}}\left[\frac{x_0^{1-\gamma}x_t^{1-\gamma}}{\sigma^2
(1-\gamma)^2 t}\right],
\end{equation}
where $\Pcal_+$ corresponds to the absorbing case and $\Pcal_-$ to
the reflecting one. In the latter case, the relationship of the
CEV and CIR processes can be exploited to sample the probability
distribution (\ref{pdfssCEV1al12}), with the same result as in
(\ref{samplingCEVagt1}). Simulation of the CEV process in this
case is shown in figure (\ref{figcev}). However, in the absorbing
case, $\Pcal_{-}$ cannot be sampled either by using the
corresponding CIR process or by using {\em negative dimensions} as
in the previous case, since for $\gamma < 1/2$ we have that
$1<1/[2(\gamma -1)]\leq 0$. Thus $\Pcal_{-}$ should be sampled
using rejection or transformation methods.

\section{Measure valued diffusions / Reaction-diffusion problems}

During the last decades much attention has been devoted to
stochastic spatial models of interacting and branching particles
like the contact model, the voter model, the normal and oriented
percolation, etc. \cite{sladeAMS,hinrichsen,carmona} Despite its
simplicity these models display interesting critical properties at
high spatial dimensions and serve as universality classes for more
complicated situations. While much analytical progress has been
reached in the study of this models, some properties of them must
be understood by using numerical methods. In this respect, several
efficient particle stochastic simulations has been proposed and
studied. An alternative approach is to study the convergence of
the particle process to a continuum measure value diffusion in
which the system is described by the (stochastic) concentration of
particles $\rho(x,t)$ at a given spatial location $x$. This
accelerates the numerical simulations and also could help to
identify the relevant dynamics at the coarse-grained dynamics.
Useful representations of these measure-valued diffusion are the
ones interpreted as solutions of a martingale problem in terms of
stochastic partial differential equations. Despite its clear
interpretation, this representation is not usually considered in
numerical simulations. The reason for that is the inability of
usual numerical methods to handle correctly the non-negativity
character and the Poissonian fluctuations of the concentration of
particles close to $\rho = 0$ \cite{moro,gaines}. In this respect,
we will see that the splitting scheme provides a very efficient
and accurate method to study these models.

\subsection{Super-Brownian motion}

The most simple and studied measured-value diffusion is the
super-Brownian motion. The super-Brownian motion arises as the
scaling limit in various critical branching systems when the
interaction between them is weak, i.e. when the system is above
some {\em critical spatial dimension} \cite{sladeAMS}. Above this
critical dimension we expect a Gaussian limit and indeed the
super-Brownian motion is the Gaussian limit of a number of models:
the voter model above 2 dimensions, the contact process above 4
dimensions, oriented percolation above 4 dimensions and
percolation over 6 dimensions.

Super-Brownian motion can be studied analytically by using the
log-Laplace transform that maps its dynamics into a non-linear
partial differential equation, a result due to Dynkin
\cite{dynkin}. More general situations or particular properties of
the sBm can only be reached through numerical simulations. To our
knowledge there is no numerical simulation of the martingale
problem of the sBm. In one dimension, the martingale problem of
the super-Brownian (sBm) motion is described by the stochastic
partial differential equation
\begin{equation}\label{eq-sbm}
d \rho(x,t) = \Delta \rho(x,t) dt + \sqrt{\sigma \rho(x,t)}
dW(x,t)
\end{equation}
where $W(x,t)$ is a Wiener sheet. It is well known that the
solutions of the sBm die in finite time almost surely. Another
interesting property is that the support of the solution, i.e. the
set for which $u(x,t) > 0$ is compact, provided that the initial
condition has a compact support.

\begin{figure}
\begin{center}
\includegraphics[width=2.5in,clip=]{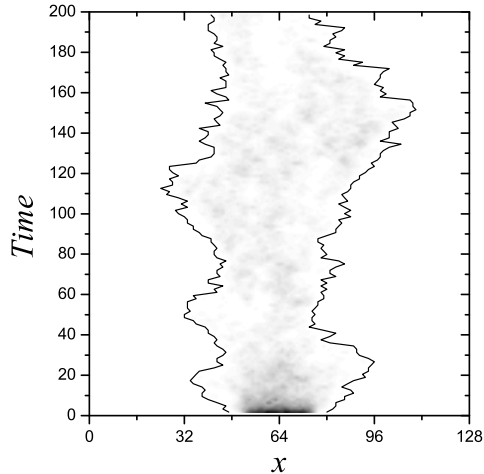}
\caption{\label{sbm1} Realization of the sBm in one dimension. The
figure shows density plots for $u(x,t)$ as a function of time with
initial condition $u(x,0) = 0.1$. The solid line depicts the
extremes of the support at each time. Parameters are $\Delta x =
1$, $\Delta t=0.1$ and $\sigma = 1.0$ in a $L=128$ lattice.}
\end{center}
\end{figure}

We use the splitting-step method for approximating the strong path
solutions of (\ref{eq-sbm}). To this end we approximate the sBm by
the super-random walk on $\Z^d$ (see \cite{carmona})
\begin{equation}\label{eq-sbm1}
du_i(t) = \Delta_i u_i(t) dt + \sqrt{\frac{\sigma u_i(t)}{ (\Delta x)^d}}
dW_i(t)
\end{equation}
where $\Delta_i$ is the discrete Laplacian operator
$$ \Delta_i u = \frac{u_{i+1}- 2 u_i + u_{i-1}}{(\Delta x)^2} $$
in $\Z^d$,
$\Delta x$ is the lattice spacing and $W_i(t)$ are independent
Wiener processes in time $t$. Equation (\ref{eq-sbm1}) represents a model of
interacting Feller diffusions. The splitting-step algorithm in
this case is based on
\begin{eqnarray}\label{eq-sbm111}
du_i^{(1)}(t)&=&\sqrt{\frac{\sigma u^{(1)}_i(t)}{(\Delta x)^d}} dW_i(t) \\
d u_i^{(2)}(t)&=&\Delta_i u_i^{(2)}
\end{eqnarray}
where the first step is integrated using the transition probability
(\ref{cdf-2}) and its sampling (\ref{rhodett}) and the
equation below can be integrated using standard schemes. In Figure
(\ref{sbm1}) we observe a simulation of the sBm in one dimension.
Our method not only provides an accurate description both in the
strong and weak sense of the sBm, but it also incorporates one of the
main properties of the sBm, namely the compact support property.

In two dimensions equation\ (\ref{eq-sbm}) is not well defined,
but still we can study equation\ (\ref{eq-sbm1}) in the lattice
$\Z^2$. As our simulations show, the compact support
properties of sBm in two dimensions are preserved and we also see
the cluster formation and their disappearance at large times. In
contrast to that, the methods of Gaines \cite{gaines} do not
provide such efficient numerical approximations.

\begin{figure}
\begin{center}
\includegraphics[width=4.5in,clip=]{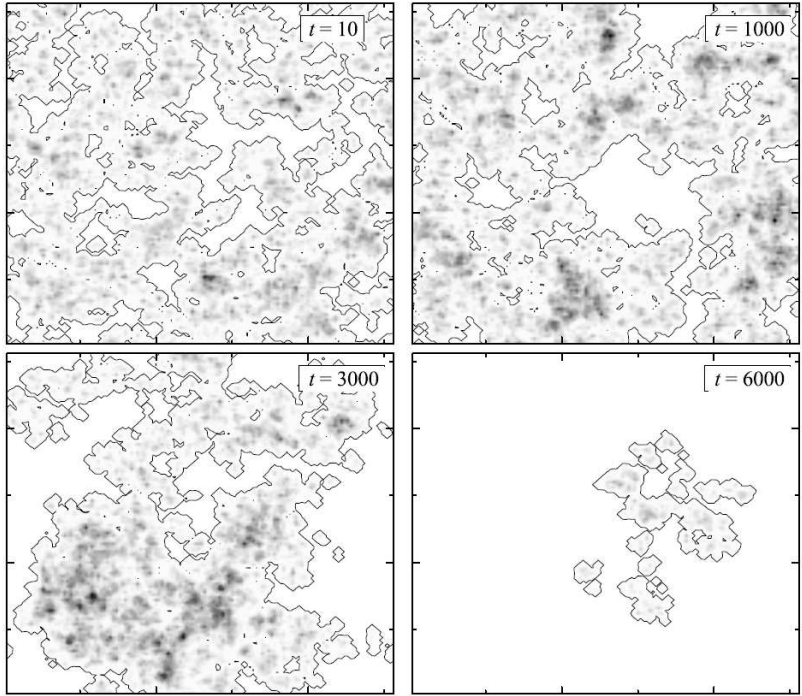}
\caption{\label{sbm2} Strong path approximation of the sBm in two
dimensions. The figures show density plots of $u(x,t)$ surrounded
by the border of the support (solid line) at different times. The
initial condition is $ u(\vec r,0) = 0.1$. Parameters are $\Delta
x = 1$, $\Delta t = 0.1$ in a $256 \times 256$ lattice.}
\end{center}
\end{figure}

\subsection{Contact process}
The contact process is a model of spreading of infection in a
lattice in which an site can be infected via contact with an
infected site in its neighborhood \cite{sladeAMS}. Infection and
recovery happens at different rates and there is a critical value
of them for which an infection started from a single infected
individual will die out in finite time. The contact process with
finite range of infection has a critical dimension $d_c=4$. At
higher dimensions the critical contact process converges to the
sBm, but not below $d_c$. At lower dimensions, it is not known
what the scaling limits should be.

However, long-range interaction with suitable scalings show that
the contact process converges to the sBm for $d \geq 2$, see
\cite{durretperkins}. In one dimension, Mueller and Tribe
\cite{mueller} showed that the rescaled density of particles for
long-ranged contact process weakly converges to the solution of
the SPDE
\begin{equation}\label{eq-contact}
d\rho(x,t) = [\Delta \rho(x,t) + \theta \rho(x,t) - \rho(x,t)^2 ]
dt + \sqrt{\rho} dW(x,t)
\end{equation}
where $W(x,t)$ denotes a space-time Wiener process.
Moreover, they showed that the above equation undergoes a phase
transition at a critical value of $\theta_c$ for which
\begin{equation}\label{def-prob}
\Pcal(u(x,t)\ \mathrm{survives}) \left\{
\begin{array}{cc}
=0 & if\ \theta < \theta_c\\
> 0 & if\ \theta > \theta_c
\end{array}\right. .
\end{equation}
The non-trivial behavior of the solution of (\ref{eq-contact}) is
believed to represent the well known universality class for
contact processes (also named directed percolation universality
class \cite{hinrichsen,moro}).

\begin{figure}
\begin{center}
\includegraphics[width=4.5in,clip=]{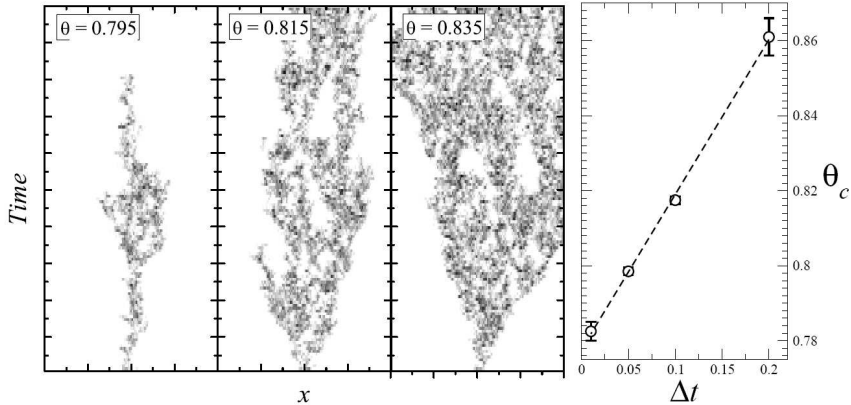}
\caption{\label{contact} (left) Strong path approximation of equation\
(\ref{eq-contact}) for $\Delta x=1$ and $L=2^{14}$ and different
values of $\theta$ below and above the phase transition. (Right)
Critical value of $\theta_c$ as defined in (\ref{def-prob}) as a
function of $\Delta t$. The straight line is a linear fit to the
data. Parameters are $\Delta x=1, L=2^{14}$.}
\end{center}
\end{figure}

As before, numerical simulations of (\ref{eq-contact}) can be now
addressed using the splitting scheme and the efficient random
number generators for the conditional probability \cite{moro}. In
particular, we discretize the spatial operators in a lattice $\Z$
[like in (\ref{eq-sbm1})] and split the dynamics as follows:
\begin{eqnarray}\label{splittingDP-1}
d\rho_i^{(1)}&=&\sqrt{\rho_i^{(1)}} dW(t)\\
d\rho_i^{(2)}&=&[\Delta_i \rho_i^{(2)} + \theta \rho_i^{(2)}
-(\rho_i^{(2)})^2]dt
\end{eqnarray}
where the last equation is numerically integrated using Euler
approximations with sufficiently small step sizes.
Results for strong approximations of $\rho_i(t)$ are shown in
Figure \ref{contact} where we can see a typical realization for
the subcritical (infection dies out), critical and super-critical
(infection spreads) for different values of $\theta$.

We can calculate the critical value of $\theta_c$ in one dimension
using our algorithm. To this end, we identify the critical point
$\theta_c$ using equation (\ref{def-prob}) and finite scaling
techniques of statistical mechanics. In particular we found that
$\theta_c = 0.777\pm0.001$ for $\Delta x=1$, and the convergence
to this value is of order one.

\section{Conclusion}
The general idea of this paper is to propose a splitting of SDE
(\ref{inicio}) into two new SDEs for which one can keep nonnegativity
during integration of both subsystems. This is achieved by either solving
one subsystem exactly in the pathwise sense or using its transition
probabilities, and solving the other (here nonrandom) subsystem by
nonnegativity-preserving numerical methods. In this way one is able to
preserve nonnegativity and a maximum of convergence order 1.0 both in weak
and strong sense. For the efficiency of our splitting algorithm, it is crucial
to find a splitting into appropriate subsystems. For this purpose, one also
tries to incorporate more complicated boundary conditions of the original
SDE (\ref{inicio}) into an explicitly solvable subsystem such that a fairly
easier numerical integration of the remaining subsystem remains to be
implemented. For example, for SDEs
\begin{eqnarray*}
dX(t) & = & [\alpha (X(t),t) + \lambda X(t)] dt + \sigma X(t) dW(t) ,
\end{eqnarray*}
one makes use of the splitting into
\begin{eqnarray*}
dX_1(t) & = & \lambda X_1(t) dt + \sigma X_1(t) dW(t) , \\
dX_2(t) & = & \alpha(X_2(t),t) dt ,
\end{eqnarray*}
where the explicit solution of the first component $X_1$ is given by
\begin{eqnarray*}
X_1(t) & = & X_1(0) \exp ((\lambda-\frac{1}{2} \sigma^2) t + \sigma W(t) )
\end{eqnarray*}
which possesses the monotone property of leaving the positive axes
$[0,+\infty)$ invariant by this type of random mapping (almost surely). This idea can
be easily extended to nonlinear systems of SDEs with its splitting into linear
and nonlinear subsystems in several dimensions. Another type of splitting is found for
nonlinear systems
\begin{eqnarray}
\label{eqSDE}
dX(t) & = & f(X(t),t) dt + \sigma (X(t),t) dW(t)
\end{eqnarray}
with continuously differentiable coefficients $\sigma$ as follows.
Rewrite this equation to as
\begin{eqnarray*}
dX(t) & = & \left[f(X(t),t) - \frac{1}{2} \sigma(X(t),t) \frac{\partial \sigma
(X(t),t)}{\partial x}  + \frac{1}{2} \sigma(X(t),t) \frac{\partial \sigma
(X(t),t)}{\partial x}  \right] dt + \\
& & \qquad + \sigma (X(t),t) dW(t)
\end{eqnarray*}
and set
\begin{eqnarray*}
\alpha(x,t) & = & f(x,t) - \frac{1}{2} \sigma (x,t) \frac{\partial \sigma (x,t)}{\partial x} ,\\
\beta (x,t) & = & \frac{1}{2} \sigma (x,t) \frac{\partial \sigma
(x,t)}{\partial x}
\end{eqnarray*}
Then the splitting-step algorithm is applied to system $(X_1,X_2)$
satisfying equations (\ref{ss-1}) and (\ref{ss-2}) with coefficients
$\alpha$ and $\beta$ as defined above. This works at least efficiently if
$\sigma(x,t)=\sigma(x)$ does not depend on $t$ and the invertible integral
$H(z) = \int^z [b(z)]^{-1} dz $ exists on the domain of definition of
the original equation (\ref{eqSDE}) for $X$. In this case one finds
$$X_1(t)=H^{-1} (W(t) + H(X_1(0)) ). $$
\par
Once an appropriate splitting is found then it is relatively easy
to implement the related numerical algorithm. The proposed splitting-step method
efficiently works since its implementation essentially relies on the well-known
variation-of-parameters formula for perturbed dynamical systems which
extends to SDEs. Recall that, by this formula, if the equation
\begin{eqnarray*}
dX(t) & = & \beta (X(t),t) dt + \sigma (X(t),t) dW(t)
\end{eqnarray*}
has known fundamental solution $\Phi=\Phi(t,X_0)$ then the exact solution
of the original equation (\ref{th11-1}) possesses the pathwise
representation
\begin{eqnarray*}
X(t+\Delta t) & = &
\Phi(t+\Delta t,X(t)) + \Phi(t+\Delta t,X(t)) \int^{t+\Delta t}_t
[\Phi(s,X(t))]^{-1} \alpha (X(s),s) \, ds\\
& \approx & \Phi(t+\Delta t,X(t)) + \alpha(X(t),t) \Delta t
\end{eqnarray*}
on each subintervals $[t,t + \Delta t] \subset [0,T]$. Thus, the motivation
of our splitting-step technique is apparent by finding $\Phi$ and
numerical integration of expressions $\int \alpha(X(s),s) \, ds$.

%\newpage
\section*{Acknowledgments}
We are grateful to C.\ Doering and P.\ Smereka for comments and
discussions. Financial support is acknowledged from the Ministerio
de Ciencia and Tecnolog\'{\i}a (Spain) through projects
FIS2004-01001 and BFM2002-04474-C02.
%The author thanks the anonymous authors whose work largely
%constitutes this sample file. He also thanks the INFO-TeX mailing
%list for the valuable indirect assistance he received.

\appendix

\section{Non-central chi-square distribution random number generation}
The splitting-step method proposed in this paper relies on the
exact numerical sampling of transition probability density for
some processes. In particular, efficient generation of random
numbers $\chi'^2_d(\lambda)$ with a non-central chi-square
distribution with $d$ degrees of freedom whose probability density
function is found in \cite{book} by noting that
\begin{equation}\label{ncchisquarepdf}
\Pcal[\chi'^2_d(\lambda) = x] = p(x;d,\lambda) =
\frac{e^{-(\lambda+x)/2}}{2}
\left(\frac{x}{\lambda}\right)^{(d-2)/4} I_{(\nu-2)/2}
(\sqrt{\lambda x}),\quad x>0 .
\end{equation}
This distribution is properly defined for any $d$ positive, and
was extended to the case $d=0$ by Siegel \cite{siegel}. Here we
will extend it to the case $d=-2,-4,\ldots$ and will show how to
sample this distribution.

To this end, we use the fact that the distribution
(\ref{ncchisquarepdf}) can be expressed also as a mixture of
central $\chi^2$ variables with Poisson weights
\begin{equation}\label{nccmixture}
p(x;d,\lambda)=\sum_{j=0}^{\infty} \frac{(\lambda/2)^j
e^{-\lambda/2}}{j!}\ p_0(x;d+2j), \quad d > 0,
\end{equation}
where $p_0(x;d)$ is the distribution of a chi-square random
variable $\chi^2_d$ with $d$ degrees of freedom. This expression
suggests a simple and efficient procedure to obtain $\chi'^2_d(\lambda)$
random variables:
\begin{enumerate}
\item Choose $K$ from a Poisson distribution with mean $\lambda/2$
so that $\Pcal[K=k] = e^{-k/2}(\lambda/2)^k/k!$ $ (k=0,1,\ldots)$.
\item Then take $\chi'^2_d(\lambda) = \chi^2_{d+2K}$, which can be
done using the any standard random number generator of the
$\chi^2_d$ distribution. %% Which citation ? \cite{}.
\end{enumerate}
In the $d=0$ case, the $\chi'^2_d(\lambda)$ distribution has a
discrete component at zero with mass $e^{\lambda/2}$ (which
represents the probability to get absorbed at zero in our
stochastic processes), see \cite{siegel}. We have
\begin{equation}\label{nccmixtured0}
p(x;0,\lambda)=\sum_{j=1}^{\infty} \frac{(\lambda/2)^j
e^{-\lambda/2}}{j!}\ p_0(x;2j) + e^{-\lambda/2}\delta(x),
\quad d = 0 .
\end{equation}
The procedure above can be modified to account for this discrete
component by taking the convention that the central $\chi^2_d$
distribution is identically zero when $d=0$. This convention can
be extended to even negative dimensions to get
\begin{equation}\label{nccmixtured1}
p(x;d,\lambda)=\sum_{j=|d|/2+1}^{\infty} \frac{(\lambda/2)^j
e^{-\lambda/2}}{j!}\ p_0(x;d+2j) + \delta(x) \sum_{j=0}^{|d|/2}
\frac{(\lambda/2)^j e^{-\lambda/2}}{j!}
\end{equation}
for $d = 0,-2,-4,\ldots$.

Summarizing, if $K$ is a Poisson random number with mean
$\lambda/2$ we have
\begin{equation}\label{samplingchi21}
\chi'^2_d(\lambda) = \chi^2_{d+2K},\qquad d>0
\end{equation}
and
\begin{equation}\label{samplingchi22}
\chi'^2_d(\lambda) = \left\{
\begin{array}{ll}0&\mathrm{if}\ d+2K \leq 0\\
\chi^2_{d+2K}&\mathrm{if}\ d+2K>0 \end{array}\right. \qquad
d=0,-2,-4,\ldots.
\end{equation}
This sampling of the probability distribution function is exact
and should be used especially when $\lambda$ is small. However,
when $\lambda$ is large, the $\chi'^2_d(\lambda)$ distribution
asymptotically converges to the Gaussian distribution and other
approximations (like the ones in \cite{book}) might be considered
to improve the speed of our algorithm.

\end{document}